\documentstyle[12pt]{article}

\catcode`\ =\active \def {\`e} \catcode`\ˆ =\active \def ˆ{\`a} \catcode`\ =\active \def {\`u}
\catcode`\Ž =\active \def Ž{\'e} \catcode`\ =\active \def {\c c} \def \parn{\par\noindent} 
\catcode`\ =\active \def {\^e} \catcode`\™ =\active \def ™{\^o} 
\catcode`\" =\active \def "{\^\i } 

\catcode`\: =\active \def :{\ifmmode\string:\else\unskip\kern 2pt\string: \ignorespaces\fi } 
\catcode`\; =\active \def ;{\ifmmode\string;\else\unskip\kern 2pt\string; \ignorespaces\fi } 

\font \db = msbm10 at 12 pt

  \def \b{\beta}   \def \g{\gamma}     \def \D{\Delta}  \def \G{\Gamma}
   \def \s{\sigma}   \def \t{\theta}  
\def \e{\varepsilon}  \def \f{\varphi}      

\def \R{\mbox{\db R}}  \def \H{\mbox{\db H}} \def \E{\mbox{\db E}} \def \P{\mbox{\db P}}

\def \S{\mbox{\db S}}      \def \Z{\mbox{\db Z}}

\def \LL{{\cal L}}  \def \M{{\cal M}}       \def \O{{\cal O}}

         \def \tt{\tilde{\theta}} 

\def \p{\partial}        \def \ii{\infty}    
     \def \ub{\underbar}

              \def \nea{\nearrow}        

   \def \ds{\displaystyle}   \def \ts{\textstyle}       
\def \rt1{\sqrt{-1}\,\,}  \def \1{^{-1}}            \def \2{^{-2}}             \def \5{{\ts {1\over 2}}}
\def \parn{\par\noindent}     \def\indf{\leavevmode\indent }

\def\cotg{{\rm cotg}\,}

\begin{document}

\newtheorem{defi}{Definition}
\newtheorem{theo}{Theorem}
\newtheorem{prop}{Proposition}  \newtheorem{propr}{Property}
\newtheorem{cor}{Corollary}
\newtheorem{lem}{Lemma}
\newtheorem{rem}{Remark}

\title{\bf Relativistic Diffusions}

\author{Jacques FRANCHI \quad and \quad Yves LE JAN}
\date { March 2004 }
\maketitle

\begin{abstract}
The purpose of this note is to introduce and study a relativistic motion whose acceleration, in proper
time, is given by a white noise. We begin with the flat case of special relativity, continue with the 
case of general relativity, and finally consider more closely the example of the Schwarzschild space.
\par  
   A detailed and completed version of this work is in progress. 
\end{abstract}

\section{Introduction} \indf 
   Several authors, mathematicians and physicists as well, have been interested for a long time in
studying a stochastic relativistic process. \par 

   In the present note, we consider diffusions, that is to say continuous Markov processes. \par 
   We start, in Section \ref{sec.relr} below, with the flat case of Minkowski space
$\,M_{1,d}\equiv\R^{1,d}\,$, and therefore with the Brownian motion of its unit pseudo-sphere, integrated
then to yield the only true relativistic diffusion, according to [D1]. We get then in a simple way its
asymptotic behavior, hopefully simplifying the point of view of [D3]. \par 
   
   In Section \ref{sec.relg} below, we present an extension of the preceding construction to the framework
of general relativity, that is to say of a generic Lorentz manifold. For that, we define first a process
at the level of pseudo-orthonormal frames, with Brownian noise only in the vertical directions, and
then see that this diffusion projects into on a diffusion on the pseudo-unit tangent bundle. The
infinitesimal generator is the sum of the vertical Laplacian with the vector field generating the geodesic
flow. \par
   In our last Section \ref{sec.S}, we deal in detail with the Schwarzschild space. This is the most
classical model for the complement in $\R^{1,3}$ of some central body, star or black hole ; see for example
[DF-C], [F-N], or  [S]. \  In this setting, the relativistic diffusion projects on a three-dimensional
diffusion. The choice of coordinates is suggested by the integration of the geodesic flow. \par 

   We show that almost surely, our diffusion hits the central body, or wanders out to infinity, both
events occuring with positive probability ; and that in the second case it goes away in some random
asymptotic direction, asymptotically with the velocity of light.

\section{A relativistic diffusion in Minkowski space} \label{sec.relr} \indf 
   Let $\,M_{1,d}\,$ be the Minkowski space of dimension $d+1$, and $\,\H^d\,$ the positive time part 
of its pseudo-unit sphere, which is a representation of the $d$-dimensional hyperbolic space.  
We identify the unitary tangent bundle $\,T^1M_{1,d}\,$ with $\,{\ds M_{1,d}\times\H^d}\,$, and we define
on it the relativistic diffusion $\,(\xi_s,B_s)\,$ by $\;{\ds\xi_s:=\xi_0+\s\int_0^sB_t\,dt}\;$, where
$\,B\,$ denotes the Brownian motion on $\,\H^d\,$ (started from some $B_0\in\H^d$), $\,\s\in\R^*\,$, and
$\,\xi_0\,$ is a fixed point in $\,M_{1,d}\,$. \par 
   Its infinitesimal generator is $\;{\ds B\,{\p\over\p\xi} +{\ts{\s^2\over 2}}\,\D^H_B}\;$, where
$\,\D^H\,$ denotes the hyperbolic Laplacian. \ It is invariant in law under any Lorentz transformation. \
In any fixed frame the time coordinate $\,\xi^0_s\,$ is strictly increasing and the velocity
$\,{\ds \Big(\,v^j_s:=\,{\dot{\xi}^j_s\over\dot{\xi}^0_s}\;\Big|\, 1\le j\le d\Big)}\,$ is bounded by 1,
the velocity of light. \par 

   From the asymptotic behavior of the hyperbolic Brownian motion, it is easy to see that there exists
almost surely some $\,\t\in\S^{d-1}\,$ such that $\,v^j_s\,$ and $\,\xi^j_s/\xi^0_s\,$ converge towards 
$\,\t^j\,$ for $\,1\le j\le d\,$ and $\,s\to +\ii\,$. \ So our diffusion process almost
surely wanders out to infinity, asymptotically in one direction and with the velocity of light. 

\section{Extension to General Relativity} \label{sec.relg} \indf 
   Let us start with a given pseudo-Riemannian manifold $\,\M\,$, whose metric tensor has signature 
$\,(1,d)\,$ at each point. Let $\,SO_{1,d}\M\,$ denote its principal bundle of pseudo-orthonormal frames,
which has its fibers modelled on the special Lorentz group.  \par 
   Let $\; H_0,..,H_d\;$ and $\,\{ V_{kl}\,|\,0\le k<l\le d\}\,$ denote the canonical horizontal and
vertical vector fields on $\,SO_{1,d}\M\,$, and denote by $\,\{ w^j_s\,|\,1\le j\le d\}\,$ $d$
independent real Wiener processes. \    The Stratonovitch stochastic differential equation 
$$ (*)\qquad du_s = H_0(u_s) ds + \s\sum_{j=1}^dV_{0j}(u_s)\circ dw^{j}_s  $$ 
has for each $\,u_0\in SO_{1,d}\M\,$ a unique solution which is a diffusion on $\,SO_{1,d}\M\,$ with 
infinitesimal generator $\;{\ds H_0+{\ts{\s^2\over 2}}\,\sum_{j=1}^dV_{0j}^2}\,$. \par 
   Let $\,\pi\,$ denote the canonical projection from $\,SO_{1,d}\M\,$ onto $\,\M\,$, and $\,\pi_1\,$
denote the canonical projection from $\,SO_{1,d}\M\,$ onto the unitary tangent bundle
$\,T^1\M\,$, which we identify with $\,SO_{1,d}\M\Big/ SO_d\,$. \quad  
   Set $\;\xi_s:= \pi(u_s)\,$, and $\; u_s=\Big(\xi_s\,; e_0(s),..,e_d(s)\Big)\,$. \par 
   Let $\,D\,$ denotes the covariant differential along the curves : in local coordinates
$\;(\xi^i,e_j^k)\,$, with $\,e_j=e_j^k{\p\over\p x^k}\,$, it writes $\;(D e_j)^k =
de_j^k+\G_{li}^ke_j^ld\xi^i\;$. \quad  The equation $\,(*)\,$ writes equivalently 
$$ \dot{\xi}_s = e_0(s) \quad ;\quad D e_0(s) = \s\sum_{j=1}^de_j(s)\circ dw^{j}_s \quad ;\quad 
D e_j(s) = \s\, e_0(s)\circ dw^{j}_s \quad \hbox{for } 1\le j\le d\; . $$ 

   The stochastic flow defined by $\,(*)\,$ commutes with the action of $\,SO_d\,$ on $\,SO_{1,d}\M\,$,
and therefore the projection $\;(\xi_s,\dot{\xi}_s) = \pi_1(u_s)\;$ is a diffusion on $\,T^1\M\,$. \ 
When $\,\M=M_{1,d}\,$, it coincides with the diffusion defined in section \ref{sec.relr} above (see Lemma 
\ref{lem.verif} below).  \par 

   So we got the following general existence result for our relativistic diffusion
$\,(\xi_s,\dot\xi_s)\,$. 
\begin{theo}\label{the.gen}\quad The $\,SO_{1,d}\M$-valued Stratonovitch stochastic differential
equation  
$$ (*)\qquad du_s = H_0(u_s) ds + \s\sum_{j=1}^dV_{0j}(u_s)\circ dw^{j}_s  $$ 
defines a diffusion $\,(\xi_s,\dot\xi_s) := \pi_1(u_s)$ on $T^1\M$, whose infinitesimal generator
is $\,{\ds\LL_0+{\ts{\s^2\over 2}}\,\D_v}\,$, $\,\LL_0\,$ denoting the generator of the geodesic flow and
$\,\D_v\,$ denoting the vertical Laplacian. 
\end{theo}
    
\section{Example of the Schwarzschild space} \label{sec.S}\indf 
   This is the most classical model for the complement of a spherical central body, star or black hole ; 
see for example [DF-C], [F-N], [S]. \par    
   We take $\;\M:= \Big\{ \xi = (t,r,\t)\in\R\times [R,+\ii [\times\S^2\Big\}\,$, where $\,R\in\R_+\,$ is
a parameter of the central body, endowed with the radial pseudo-metric : 
$$ (1-{\ts{R\over r}})\, dt^2 - (1-{\ts{R\over r}})\1 dr^2 - r^2 |d\t|^2\, . $$ 
The coordinate $\,t\,$ represents the time (multiplied by the velocity of
light $\,c\,$), and $r$ the distance from the origin. 

   In spherical coordinates $\,\t =(\f,\psi )\in [0,\pi ]\times (\R/2\pi\Z )$, we have 
\mbox{$\,|d\t|^2\! =d\f^2+\sin^2\!\f d\psi^2$,} \ and the non-vanishing Christoffel symbols are :
$$ \G_{rt}^t = -\G_{rr}^r = {R\over 2r(r-R)}\; ;\; \G_{tt}^r = {R(r-R)\over 2r^3}\; ;\; 
\G_{\f\f}^r = R-r\; ;\; \G_{\psi\psi}^r = (R-r)\sin^2\f\; ;\; $$ 
$$ \G_{r\f}^\f = \G_{r\psi}^\psi = r\1\; ;\; \G_{\psi\psi}^\f = -\sin\f\cos\f\; ;\;
\G_{\f\psi}^\psi = \cotg\f\; . $$ 
The Ricci tensor vanishes, the space $\,\M\,$ being empty. \ The limiting case $\,R=0\,$ is the flat case
of special relativity, considered in section \ref{sec.relr}. \  There is no other radial pseudo-metric 
in $\,\M\,$ which satisfies these constraints. 

\begin{lem}\label{lem.verif}\quad The case $R=0$ of the Schwarzschild space is actually the case of
special relativity studied in Section \ref{sec.relr}.  
\end{lem} 
\ub{Proof}\quad Using Theorem \ref{the.gen}, we just have to to check that the vertical
Laplacian is actually the hyperbolic Laplacian. Now, owing to the triviality of the tangent bundle, this
vertical Laplacian is here merely the restriction to the unit pseudo-sphere $\,\H^3\,$ of the Laplacian on
the Lorentz group, induced by the ``boost'' transformations, which is indeed the Laplacian of the
hyperboloid $\,\H^3\,$. $\;\diamond$ \par
\smallskip 

   If $\;(\xi_s,\dot{\xi}_s) = \pi_1(u_s)\;$ is the relativistic diffusion of Section \ref{sec.relg}, we
easily see that \parn 
$\Big( r_s,T_s:=\dot{r}_s,U_s:=|\dot\t_s|\Big)\,$ is an autonomous diffusion, \  with lifetime \parn 
$D:=\min\{ s> 0\,|\, r_s = R\}\,$, and infinitesimal generator 
$$ T{\p\over\p r} + \s^2TU{\p^2\over\p T\p U} + {\ts{\s^2\over 2}}\, (T^2+1-{\ts{R\over r}}){\p^2\over\p
T^2}  + {\ts{\s^2\over 2}}\, (U^2+{\ts{1\over r^2}}){\p^2\over\p U^2} \, + \hskip 30mm $$ \vspace{-4mm} 
$$ \hskip 30mm  +\, \Big({\ts{3\s^2\over 2}}\,T + (r-{\ts{3\over 2}}R)\, U^2-{\ts{R\over 2r^2}}\Big)
{\p\over\p T} + \Big( {\ts{3\s^2\over 2}}\,U-{\ts{2\,TU\over r}}+{\ts{\s^2\over 2\,r^2U}}\Big)
{\p\over\p U}\; . $$
Note that the unit pseudo-norm relation writes \quad 
${\ds (1-{\ts{R\over r_s}})(\dot t_s)^2-(1-{\ts{R\over r_s}})\1T_s^2 = r_s^2U_s^2 +1}\,$. \par 

   If we set (according for example to ([F-N], 4.4) in the deterministic case) : 
$$ a_s := (1-R/r_s)\, \dot t_s \; , \quad \hbox{ and }\quad b_s:= r_s^2\, U_s\; , $$ 
we get \ $\;T^2=a^2-(1-R/r)(1+b^2/r^2)$ , \ and the following : 
\begin{prop}\label{pr.diff1}\quad The process $\,(r_s,a_s,b_s,T_s)\;$ is a degenerate diffusion,
with lifetime $\; D=\min\{ s> 0\,|\, r_s = R\}$, and infinitesimal generator 
$$ \LL := T{\p\over\p r} + {\s^2\over 2}\,\Big( a^2-1+{\ts{R\over r}}\Big)\,{\p^2\over\p a^2} +
{\s^2\over 2}\, (b^2 + r^2) \,{\p^2\over\p b^2} + \s^2ab\,{\p^2\over\p a\p b} \hskip 30mm $$ 
$$ + \,{{3\s^2\over 2}}\, a \,{\p\over\p a} + \Big({{3\s^2\over 2}}\, b + {\s^2\, r^2\over 2 b}\Big)\,
{\p\over\p b} + \s^2aT\,{\p^2\over\p a\p T} + \s^2bT\,{\p^2\over\p b\p T}  \hskip 2mm $$
\vspace{-4mm} 
$$ \hskip 30mm  + \,{\s^2\over 2}\,\Big( T^2+1-{\ts{R\over r}}\Big)\,{\p^2\over\p T^2} +
\Big({{3\,\s^2\over 2}}\, T + (r-{\ts{3\over 2}}R)\,{b^2\over r^4} - {R\over 2r^2}\,\Big){\p\over\p T} 
\; . $$ 
Equivalently, we have the following system of stochastic differential equations : 
$$ dr_s = T_s\, ds\; ,  \quad dT_s =  dM^1_s + {\ts{3\,\s^2\over 2}}\, T_s\, ds + 
(r_s- {\ts{3\over 2}}R)\,{b_s^2\over r_s^4}\, ds - {R\over 2r_s^2}\, ds \; , $$ 
$$ d a_s = dM^a_s + {\ts{3\,\s^2\over 2}}\, a_s\, ds  \; , \quad  
d b_s = dM_s^b + {\ts{3\,\s^2\over 2}}\, b_s\, ds + {\s^2\, r_s^2\over 2\, b_s}\, ds \; , $$ 
with quadratic covariation matrix of the local martingale $\, (M^a, M^b, M^1)\,$ given by 
$$ K'_s := \s^2\pmatrix { a_s^2-1+{R\over r_s} & a_s\,b_s & a_s\,T_s \cr 
a_s\,b_s & b_s^2+r_s^2 & b_s\,T_s \cr a_s\,T_s & b_s\,T_s & 
T^2_s+1-{\ts{R\over r_s}}\cr } . $$ 
\end{prop}
\begin{cor}\label{cor.diff1}\quad The process $\,(r_s,b_s,T_s)\;$ is a diffusion,
with lifetime $\,D\,$ and infinitesimal generator 
$$ \LL' := T{\p\over\p r} + {\s^2\over 2}\, (b^2 + r^2) \,{\p^2\over\p b^2} + {\s^2\over 2b}\,(3b^2 + r^2)
\,{\p\over\p b} + \s^2bT\,{\p^2\over\p b\p T}  \hskip 20mm $$
\vspace{-4mm} 
$$ \hskip 30mm  + \,{\s^2\over 2}\,\Big( T^2+1-{\ts{R\over r}}\Big)\,{\p^2\over\p T^2} +
\Big({{3\,\s^2\over 2}}\, T + (r-{\ts{3\over 2}}R)\,{b^2\over r^4} - {R\over 2r^2}\,\Big){\p\over\p T} 
\; . $$
\end{cor}

   We have the following result on the behavior of the coordinate $\,a_s\,$. 
\begin{lem}\label{lem.cv}\quad There exist a standard real Brownian motion $\,w_s\,$, and a real 
process $\,\eta_s\,$, almost surely converging in $\R$ as $s\nea D\,$, such that 
$\; |a_s| = \exp (\s^2\, s+\s\, w_s+\eta_s)\;$ for all $\, s\in[0,D[\,$. \ In particular $a_s$ almost
surely cannot vanish, so that we can assume $a_s>0\,$ for all $\, s\in[0,D[\,$.  
\end{lem} 
\ub{Proof}\quad Proposition \ref{pr.diff1} above shows that \ 
$\;(a_s^2-1)\,\s^2\,ds\le \langle dM^a_s\rangle\le a_s^2\,\s^2\,ds\,$, for $\,0\le s <D\,$. \parn  
So that we have almost surely (as $s\to\ii\,$, when $D=\ii$) : 
$$ \log |a_s| - \log |a_0| \,=\; 3\s^2 s/2-\5\int_0^sa_t\2\langle dM_t^a\rangle + \int_0^sa_t\1\, dM_t^a
\;\ge\; \s^2 s + \int_0^sa_t\1\, dM_t^a $$ 
$$ = \s^2 s + o\Big( \int_0^sa_t\2\langle dM_t^a\rangle\Big) = \s^2 s + o(s) \; . $$ 
Since $\,(1-{\ts{R\over r_s}})\le a_s^2\,$, this implies 
$\;{\ds\int_0^D (1-{\ts{R\over r_s}})\, a_s\2\, ds\; <\ii}\;$ almost surely. \par 
   For some real Brownian motion $\,w\,$, the process $\,\eta\,$ being defined by the formula in the 
statement, we have : 
$$ d\eta_s = d(\log |a_s|) -\s^2\,ds-\s\,dw_s = \5\,(1-R/r_s)\, a_s\2\,\s^2\, ds +
\Big(\sqrt{1-(1-R/r_s)a_s\2 }-1\Big)\,\s\, dw_s \; , $$ 
and then for any $\,s<D\,$ : 
$$ \eta_s = \eta_0 + {\ts{\s^2\over 2}}\int_0^s (1-R/r_t)\, a_t\2\, dt - \s\int_0^s {(1-R/r_t)\,
a_t\2\over 1+ \sqrt{1-(1-R/r_t)\, a_t\2 }}\, dw_t \; , $$ 
which converges almost surely in $\R$ as $s\nea D\,$, since almost surely for all $\,s\in ]0,D[\,$ : 
$$ \Big\langle\int_0^s {(1-R/r_t)\, a_t\2\over 1+ \sqrt{1-(1-R/r_t)\, a_t\2 }}\, dw_t\Big\rangle \le 
\int_0^s \Big( (1-{\ts{R\over r_t}})\, a_t\2\Big)^2\, dt < \int_0^D (1-{\ts{R\over r_t}})\, a_t\2\,
ds < \ii\; . \;\;\diamond $$

   The deterministic case $\,\s=0\,$ corresponds to the case of the geodesic flow, so we recover that 
the functionals $\,a\,$ and $\,b\,$ are constants of motions.  \par 
   There are five types of timelike geodesics : \parn 
- running from $R$ to $+\ii$, or in the opposite direction ; \parn   
- running from $R$ to $R$ ; \parn  
- running from $+\ii$ to $+\ii$ ; \parn 
- running from $R$ to some $R_1$ or from $R_1$ to $+\ii$, or idem in the opposite direction ; \parn  
- bounded geodesics. \par  

   In the stochastic case $\,\s\not= 0\,$ we have the following. 

\begin{theo}\label{the.asymp}\quad 1) \ For any initial condition, the radial process $\,(r_s)\,$ almost
surely reaches $R$ within a finite time $D$ or goes to $+\ii$ as $s\to +\ii$ (equivalently as $\, t(s)\to
+\ii$).
\parn 
2) \ Both events in 1) above occur with positive probability, from any initial condition. \parn 
3) \ Conditionally to the event $\,\{ D=\ii\}\,$ of non-reaching the
central body, the Schwarzschild diffusion $\,(\xi_s,\dot\xi_s)\,$ goes almost surely to infinity in some 
random asymptotic direction of $\,\R^3$, asymptotically with the velocity of light. 
\end{theo} 

   In the proof of this theorem, below, we shall use the following very simple lemma.
\begin{lem}\label{lem.mart}\quad Let $\,M_\cdot\,$ be a continuous local martingale, and $\,A_\cdot\,$ 
a process such that \parn 
$\liminf_{s\to\ii}\limits\; A_s/\langle M\rangle_s >0\,$ almost surely on $\,\{ \langle M\rangle_\ii
=\ii\}$. \ Then $\;\lim_{s\to\ii}\limits\, (M_s+A_s) = +\ii\;$ almost surely on $\,\{ \langle
M\rangle_\ii =\ii\}$. 
\end{lem}
\ub{Proof} \quad Writing $\;M_s = W(\langle M\rangle_s)\;$, for some real Brownian motion $\,W$, we find
almost surely some $\,\e >0\,$ and some $\,s_0\ge 0\,$ such that $\;A_s\ge 2\e\,\langle M\rangle_s\;$ and
$\;|M_s|\le\e\,\langle M\rangle_s\;$ for $\,s\ge s_0\,$. Whence $\;M_s+A_s \ge\e\;\langle M\rangle_s\;$
for $\,s\ge s_0\,$. $\;\diamond $ 
\medskip  
 
\ub{Proof of Theorem \ref{the.asymp}} \quad  We prove successively the 3 assertions of the statement. \par 
\smallskip 
   1) \ {\it Almost sure convergence on $\,\{ D=\ii\}$ of $\,\,r_s\,$ to $\,\ii$.}
\par\smallskip 
   This delicate proof will be split into six parts. \par     
   Let us denote by $\,A\,$ the set of paths with infinite lifetime $\,D\,$ such that the radius $\,r_s\,$
remains bounded. \ We have to show that it is negligible for any initial condition 
$\,x=(r,b,T)=(r_0,b_0,T_0)\,$ belonging to the state space $\,[R,\ii[\times\R_+\times\R\,$. 
\par\smallskip 

   The cylinder $\;\{ r=3R/2\}\;$ plays a remarquable r\^ole in Schwarzschild geometry. In particular, it
contains light lines. We see in the following first part of proof that we have to deal with this
cylinder.  \par\smallskip 

   $(i)$ \ {\it $\,r_s\,$ must converge to $\,3R/2\,$, almost surely on $\,A\,$.} \par \smallskip 
   Let us apply It\^o's formula to $\; Y_s:= (1-{3R\over 2r_s})\,{T_s\over a_s}\,$ : 
$$ Y_s = M_s + {\ts{3R\over 2}}\int_0^s {T^2_t\over a_t\,r^2_t}\, dt + \int_0^s (1-{\ts{3R\over 2r_t}})^2
{b^2_t\over a_t\,r^3_t}\, dt  - \,\s^2\int_0^s (1-{\ts{R\over r_t}})\,{Y_t\over a_t^2}\, dt - 
\int_0^s (1-{\ts{3R\over 2r_t}})\,{R\over 2a_t\,r_t^2}\, dt \; , $$  
with some local martingale $\,M\,$ having quadratic variation : 
$$ \langle dM_s\rangle\, =\, (1-{\ts{3R\over 2r_s}})^2(1-{\ts{R\over r_s}})\Big( 1-{T_s^2\over
a_s^2}\Big){\s^2\over a_s^2}\; ds \;\le\; \s^2\,a_s\2\,ds\; . $$ 

   Now observe from the unit pseudo-norm relation (Property \ref{pr.diff1}, 1) that $\,{T_s/a_s}\,$
and $\,Y_s\,$ are bounded by 1. Hence Lemma \ref{lem.cv} implies that the last two terms in the expression
of $\,Y_s\,$ above have almost surely finite limits  as $\,s\to\ii\,$. Idem for $\,\langle M_s\rangle\,$,
and then for $\,M_s\,$. Moreover the two remaining bounded variation terms in the expression of $\,Y_s\,$
above increase. As a consequence, we get that $\,Y_s\,$, $\,{\ds\int_0^s {T^2_t\over a_t\,r^2_t}\, dt}\,$,
and $\,{\ds\int_0^s(1-{\ts{3R\over 2r_t}})^2\, {b^2_t\over a_t\,r^3_t}\, dt}\,$ converge almost surely in
$\R$ as $\,s\to\ii\,$. So does also $\,{\ds\int_0^s {dt\over a_t\,r^2_t}}\,$. \par 

   Now using that $\; {\ds {a\over r^2}\le \Big( a+{R\over r}\,\Big({b^2\over a\,r^2}+{1\over
a}\Big)\Big)\, r\2 = {T^2\over a\,r^2}+{b^2\over a\,r^4}+{1\over a\,r^2} }\;$, 
we deduce that almost surely 
$$ \int_0^\ii(1-{\ts{3R\over 2r_t}})^2\, \Big|{d\over dt}\,(1/r_t)\Big|\, dt =
\int_0^\ii(1-{\ts{3R\over 2r_t}})^2\, {|T_t|\over r^2_t}\, dt \le \int_0^\ii  (1-{\ts{3R\over 2r_t}})^2\,
{a_t\over r^2_t}\, dt <\ii\; . $$ 
This implies the almost sure convergence of $\;{ \Big( 1-{3R\over 2r_s}+{3R^2\over 4r_s^2}\Big)\Big/
r_s }\;$, and therefore of $\,(1/r_s)\,$. \par 
  Since $\,\lim_{s\to\ii}\limits(1/r_s)\,$ cannot be 0 on $A\,$, we have necessarily 
$\,\lim_{s\to\ii}\limits r_s = 3R/2\,$ almost surely on $A\,$, from the convergence of 
$\,{\ds\int_0^\ii  (1-{\ts{3R\over 2r_t}})^2\,{a_t\over r^2_t}\, dt }\;$. \par \bigskip

   $(ii)$ \ {\it $\;{b_s/a_s}\,$ converges to $\,3R\sqrt{3}/2\,$, and $\,T_s/b_s\,$ goes to 0, almost
surely on $A\,$.} \par \smallskip 
   Indeed, It™'s formula gives (for some real Brownian motion $w$) 
$$ {b^2_s\over a^2_s} = {b^2_0\over a^2_0}+ 2\s\int_0^s {b_s\over a_s}\,\sqrt{{r_s^2\over
b_s^2}-{1-{R\over r_s}\over a_s^2}}\, dw_s +2\s^2\int_0^s{r_s^2\over a_s^2}\,ds - 3\s^2\int_0^s (1-{R\over
r_s})\,{b^2_s\over a^4_s}\,ds \; . $$  
Since by the unit pseudo-norm relation we have $\;{b^2_s\over
a^2_s} < r_s^2/(1-{R\over r_s})\;$,  whence $\;{b^2_s/a^2_s}\;$ bounded on $A$, the above formula and
Lemma \ref{lem.cv} imply the almost sure convergence of $\;{b^2_s/a^2_s}\,$ on $A\,$. Indeed the bounded
variation terms converge, and as $\,{b^2_s/a^2_s}\,$ is positive, the martingale part has to converge also.
Using the unit pseudo-norm relation again, we deduce that 
$\;{\ds {T^2_s\over a^2_s} = 1 - (1-{\ts{R\over r_s}})\Big({b^2_s\over a^2_sr_s^2}+{1\over a^2_s}\Big)}\,$
has also to converge, necessarily to $0$, since otherwise we would have an infinite limit for $T_s\,$,
which is clearly impossible on $A$. The value of the limit of $\,b_s/a_s\,$ follows now directly from
this and from $(i)$.  
\par \medskip 

   $(iii)$ \ {\it We have almost surely on $A\,$ : $\;{\ds\int_0^\ii (r_t-{\ts{3R\over 2}})^2\,
b^{2}_t\,dt <\ii }\,$, and \ $\;{\ds\int_0^\ii T_t^2\,dt<\ii }\,$. 
} \par \smallskip 

   Let us write It\^o's formula for  
$\; Z_s:= (r_s-{\ts{3R\over 2}})\,T_s = \5\,{d\over d_s}\,(r_s-{\ts{3R\over 2}})^2\,$ : 
$$ Z_s = Z_0 + M_s + {\ts{3\s^2\over 2}}\,(r_s-r_0) +
\int_0^s T_t^2\,dt + \int_0^s (r_t-{\ts{3R\over 2}})^2 b^{2}_t\,{dt\over r^4_t} 
- {\ts{R\over 2}}\int_0^s (r_t-{\ts{3R\over 2}})\,{dt\over r^2_t} \, , $$  
where $\,M_\cdot\,$ is a local martingale having quadratic variation given by : 
$$\; \langle M\rangle_s = \s^2\int_0^s (r_t-{\ts{3R\over 2}})^2\, \Big( 1-{\ts{R\over
r_t}} + T_t^2\Big)\,dt \, .  $$ 

   Note that if $\,\langle M\rangle_\ii = \ii\,$, then by $\,(ii)\,$ above 
$\;{\ds\lim_{s\to\ii}\limits\; \int_0^s (r_t-{\ts{3R\over 2}})^2\, b^{2}_t\,{dt\over r^4_t}\Big/\langle
M\rangle_s =\ii}\,$. \parn  
Note moreover that in this case 
$\;{\ds \int_0^s |r_t-{\ts{3R\over 2}}|\,{dt\over r^2_t}\,\le  
\sqrt{\int_0^s (r_t-{\ts{3R\over 2}})^2\, b^{2}_t\,{dt\over r^4_t}}\times\sqrt{\int_0^s{dt\over b^2_t}}
}\;$ is also negligible with respect to $\;{\ds\int_0^s (r_t-{\ts{3R\over 2}})^2\,
b^{2}_t\,{dt\over r^4_t} }\;$. \par 
   On the other hand, we must have 
$\;{\ds\liminf_{s\to\ii}\limits\; |Z_s|= 0}\;$ on $A\,$. \par

   Therefore we deduce from Lemma \ref{lem.mart} that necessarily $\,\langle M\rangle_\ii < \ii\,$, and
then that $\,M_s\,$ has to converge, almost surely on $A\,$. \par 
   Using again that $\;{\ds\liminf_{s\to\ii}\limits\; |Z_s|= 0}\,$, we deduce the almost sure
boundedness and convergence on $A\,$ of $\;{\ds\int_0^\ii T_t^2\,dt}\;$ and of $\;{\ds\int_0^\ii
(r_t-{\ts{3R\over 2}})^2\, b^{2}_t\,{dt\over r^4_t} }\;$. 
\par\medskip 

   $(iv)$ \ {\it $(r_s-{\ts{3R\over 2}})^2 b_s\,$ and  $\,{T^2_s/b_s}\;$ go to 0 as $s\to\ii\,$, 
almost surely on $\,A\,$.} \par\smallskip  
   Indeed, on one hand we deduce from $(iii)$ that (for some real Brownian motion $\,W_\cdot$)
$$ (r_s-{\ts{3R\over 2}})^2 b_s= \s W\Big[\int_0^s (r_t-{\ts{3R\over 2}})^4(b_t^2+r_t^2)dt\Big] + 
2\!\int_0^s (r_t-{\ts{3R\over 2}})T_tb_tdt +{\ts{\s^2\over 2}}\!\int_0^s (r_t-{\ts{3R\over 2}})^2
(3b_t+{\ts{r_t\over b_t}})dt $$ 
has to converge almost surely on $\,A\,$ as $s\to\ii\,$, necessarily to 0 since it is integrable
with respect to $s\,$.  \par 
  On the other hand we have for some real Brownian motion $\,W'_\cdot$, by It™ formula : 
\parn 
\vbox{
$$ {T^2_s\over b_s} = {T^2_0\over b_0} + \s\, W'\Big[\int_0^s \Big({T^4_t\over b^2_t}+r_t^2\,{T^4_t\over
b_t^4}+4(1-{\ts{R\over r_t}}){T^2_t\over b^2_t}\Big) dt\Big] + {\ts{\s^2\over 2}}\int_0^s {T^2_t\over
b_t}\,dt + 2\int_0^s (r_t-{\ts{3R\over 2}})\,T_t\,b_t\,{dt\over r^4_t} \hskip 6mm  $$ 
$$ \hskip 50mm   - \int_0^s {R\,T_t\over r^2_t\,b_t}\,dt + \s^2\!\int_0^s(1-{\ts{R\over r_t}})\,{dt\over
b_t} +  {\ts{\s^2\over 2}}\int_0^s {r_t^2\,T_t^2\over b_t^3}\,dt \; . $$ }
Recall from $(i)$ that $\,{T_t\over b_t}\to 0\,$ and that $\,b_t\sim {3R\sqrt{3}\over 2}\,a_t\,$. 
Thus using $(iii)$ we see easily that all integrals in the above formula converge. Hence we deduce the
almost sure convergence of $\,{\ds s\mapsto {T^2_s/b_s}}\,$ on $\,A\,$, necessarily to $0\,$, since it
is integrable.  \par\medskip 

   $(v)$ \ {\it It is sufficient to show that $\;{\ds\int_0^\ii |r_t-{\ts{3R\over 2}}|\, |T_t|\,
b^{2}_t\,dt <\ii}\;$, and that $\;{\ds \int_0^\ii T_t^4\,dt <\ii}\;$, almost surely on $\,A\,$. }
\par\smallskip  
   Indeed, assuming that these 2 integrals are finite, It™'s formula shows that we have for some real
Brownian motion $\,W''_\cdot$ : \parn 
\vbox{
$$ T_s^2 = T_0^2 + 2\s\, W''\Big[\int_0^s (T_t^2+1-{\ts{R\over r_t}})\,T_t^2\,dt\Big] +
4\s^2\!\int_0^sT_t^2\,dt + 2\!\int_0^s (r_t-{\ts{3R\over 2}})\,T_t\,b^{2}_t\,{dt\over r^4_t}\hskip 10mm $$
\vspace{-3mm} 
$$ \hskip 50mm + \,\s^2\!\int_0^s(1-{\ts{R\over r_t}})\,dt - R\!\int_0^s{T_t\over r_t^2}\,dt\;  $$ }
$$ = \g_s + \s^2\int_0^s(1-{\ts{R\over r_t}})\,dt - R\!\int_0^s{T_t\over r_t^2}\,dt 
= \g_s + \int_0^s\Big[ {\ts{1\over 3}}+ {\ts{2\over 3r_t}}(r_t-{\ts{3R\over 2}})\Big] dt 
+ {R\over r_s}-{R\over r_0}\, = \g'_s + s/3\, , $$  
where $\,\g_\cdot\;,\;\g'_\cdot\,$ are bounded converging processes on $\,A\,$. \ Whence 
$\;\lim_{s\to\ii}\limits T^2_s=\ii\;$ almost surely on $A\,$, which with $(iii)$ above implies 
that $\,A\,$ must be negligible.  \par\medskip  

   $(vi)$ \ {\it End of the proof of the convergence of $\,\,r_s\,$ to $\,\ii$ on $\,\{ D=\ii\}\,$. }
\par\smallskip 
   By Schwarz inequality, the first bound in $(v)$ above will follow from 
$\;{\ds\int_0^s T_t^2b_t\,dt<\ii}\;$ and from 
$\;{\ds\int_0^s (r_t-{\ts{3R\over 2}})^2 b^{3}_t\,dt <\ii}\;$. Now these two terms appear in the It™
expression for \ $\;Z_s^1:= (r_s-{\ts{3R\over 2}})\,T_s\,b_s\;$ :  
$$ Z_s^1 = Z_0^1 + M_s^1 + {\ts{\s^2\over 2}}\int_0^s \!\Big[ 8+{r_t^2\over b_t^2}\Big] Z_t^1 dt +
\int_0^s T_t^2b_t\,dt + \int_0^s (r_t-{\ts{3R\over 2}})^2 b^{3}_t\,{dt\over r^4_t} 
- {\ts{R\over 2}}\!\int_0^s (r_t-{\ts{3R\over 2}})b_t\,{dt\over r^2_t} \;, $$ 
with a local martingale $\,M^1_\cdot\,$ having quadratic variation :  
$$ \langle M^1\rangle_s = \s^2\int_0^s (r_t-{\ts{3R\over 2}})^2\,b^{2}_t\times \Big( 1-{\ts{R\over r_t}}
+[4+ r_t^2b_t\2]\,T_t^2\Big)\,dt\; .  $$ 

      Note that by Schwarz inequality, $(iii)$ above implies that $\;{\ds\int_0^\ii
|Z^1_t|\,dt <\ii }\;$, and then that $\;{\ds \int_0^s \Big[ 8+{r_t^2\over b_t^2}\Big] Z_t^1 dt}\;$ is 
bounded and converges, almost surely on $\,A\,$, as $s\to\ii\,$. \par 
   Using the first assertion of $(iv)$, observe that 
$\;{\ds\lim_{s\to\ii}\limits\; {{\ds\int_0^s} (r_t-{\ts{3R\over 2}})^2\, b^{3}_t\,{dt\over r^4_t}
+{\ds\int_0^s} T_t^2b_t\,dt\over\langle M^1\rangle_s} =\ii}\;$ if $\,\langle M^1\rangle_\ii = \ii\,$.
\quad  Note moreover that in this case 
$$ \Big|\int_0^s (r_t-{\ts{3R\over 2}})\,b_t\,{dt\over r^2_t}\,\Big| \le 
\sqrt{\int_0^s (r_t-{\ts{3R\over 2}})^2\, b^{3}_t\,{dt\over r^4_t}}\times\sqrt{\int_0^s{dt\over b_t}} $$ 
is also negligible with respect to $\;{\ds\int_0^s (r_t-{\ts{3R\over 2}})^2\,
b^{3}_t\,{dt\over r^4_t}+{\ds\int_0^s} T_t^2\,b_t\,dt }\;$. \par 
   Therefore we deduce from Lemma \ref{lem.mart} and from the integrability of $\,t\mapsto|Z^1_t|\,$, that
necessarily $\,\langle M^1\rangle_\ii < \ii\,$, and then that $\,M^1_s\,$ has to converge, almost surely
on $A\,$. \par 
   Hence $\,Z^1_\cdot\,$ must have a limit almost surely on $A\,$, which must be 0, owing to the 
integrability of $\,Z^1_\cdot\,$. This forces clearly  
$\;{\ds\int_0^\ii (r_t-{\ts{3R\over 2}})^2\, b^{3}_t\,{dt\over r^4_t}+{\ds\int_0^\ii} T_t^2\,b_t\,dt }\;$
to be finite, almost surely on $A\,$, showing the first bound in $(v)$ above. \par \smallskip 

   Finally, the integrability of $\,T_t^2\,b_t\,\,$ and the second convergence of $(iv)$ imply 
the second bound in $(v)$ above : $\;{\ds \int_0^\ii T^4_t\,dt <\ii}\;$ almost surely on $\,A\,$. \par 
   This concludes the proof of the first assertion in Theorem \ref{the.asymp}.
\par\bigskip 

   2) \ {\it $r_s\to R\,$ and $\,r_s\to\ii\,$ occur both with positive probability, from any initial
condition. }  \par\smallskip    

   Let us use the support theorem of Stroock and Varadhan (see for example ([I-W], Theorem VI.8.1)) to show
that the diffusion $\,(r_\cdot,b_\cdot,T_\cdot)\,$ of Corollary \ref{cor.diff1} is irreducible.  
Since we can decompose further the equations given in Property \ref{pr.diff1} for
$\,(r_\cdot,b_\cdot,T_\cdot)\,$, using a standard Brownian motion $\,(w_\cdot,\b_\cdot,\g_\cdot)\in\R^3\,$,
as follows :
$$ dr_s = T_s\, ds\;,\quad d b_s = b_s\,dw_s + r_s\,d\b_s + {\ts{3\,\s^2\over 2}}\, b_s\,ds + 
{\s^2\, r_s^2\over 2\, b_s}\,ds\; , $$ 
$$ dT_s =  T_s\,dw_s + \sqrt{1-{\ts{R\over r_s}}}\,d\g_s + {\ts{3\,\s^2\over 2}}\, T_s\, ds + 
(r_s-{\ts{3\over 2}}R)\,{b_s^2\over r_s^4}\, ds - {R\over 2r_s^2}\, ds \; , $$ 
we see that trajectories moving the coordinate $\,b_\cdot\,$ without changing the others, and trajectories
moving the coordinate $\,T_\cdot\,$ without changing the others, belong to the support of 
$\,(r_\cdot,b_\cdot,T_\cdot)\,$. Moreover we see that there are
timelike geodesics,  and then trajectories in the support, which link $\,r\,$ to $\,r'\,$, and then
considering the velocities also, which link say  $\,(r,b'',T'')\,$ to $\,(r',b'',T''')\,$. \ So, for given
$\,(r,b,T)\,$ and $\,(r',b',T')\,$ in the state space, we can, within the support of
$\,(r_\cdot,b_\cdot,T_\cdot)\,$, move
$\,(r,b,T)\,$ to $\,(r,b'',T'')\,$, then $\,(r,b'',T'')\,$ to $\,(r',b'',T''')\,$, and finally move 
$\,(r',b'',T''')\,$ to $\,(r',b',T')\,$, thereby showing the irreducibility of
$\,(r_\cdot,b_\cdot,T_\cdot)\,$. 
\par 

   This implies that it is enough to show that for large enough $\,r_0,T_0\,$, the convergence to
$\ii$ occurs with probability $\ge 1/2\,$, and that for $\,r_0\,$ close enough from $R$ and $\,T_0\,$
negative enough, the convergence to $R$ occurs with probability $\ge 1/2\,$ as well. Now this can be done
by a classical supermartingale argument using the process $\,1/|T_s|\,$, stopped at some hitting time.
Indeed we see from Property \ref{pr.diff1}  that 
$$ {1\over |T_s|} + \int_0^s \Big( {\ts{\s^2\over 2}}\,T_t^2 - \s^2\,(1-{\ts{R\over r_s}}) 
- {R\,T_t\over 2\,r_t^2}\, + (2r_t-3R)\,{b_t^2\,T_t\over r_t^4}\,\Big)\, {dt\over |T_t|^3} $$  
is a local martingale. \par 
   Take first $\;r_0\ge 3R/2\,$, $\,T_0\ge 4+{4\over R\s^2}\,$, and $\;\tau :=\inf\Big\{ s>0\,\Big|\,
T_s=2+{2\over R\s^2}\,\Big\}$ : \ $r_s$ increases on $\,\{0\le s<\tau\}\,$ and then we see that 
$\,1/|T_{s\wedge\tau}|\,$ is a supermartingale, which implies that \parn 
$(2+{2\over R\s^2})\1\P (\tau <\ii )\le \liminf_{s\to\ii}\limits\, \E ({1\over |T_{s\wedge\tau}|}
1_{\{\tau <\ii\}}) \le \liminf_{s\to\ii}\limits\, 
\E ({1\over |T_{s\wedge\tau}|}) \le \E ({1\over T_{0}}) \le (4+{4\over R\s^2})\1\,$, \parn 
and then that $\;\P (\lim_{s\to\ii}\limits r_s = +\ii ) \ge \P (\tau =\ii ) \ge 1/2\,$.  \par 

   Conversely take $\;r_0\le 3R/2\,$, $\,T_0\le -{2}\,$, and $\;\tau':=\inf\Big\{ s>0\,\Big|\,
T_s=-\sqrt{2}\,\Big\}$ : \ $r_s$ decreases on $\,\{0\le s<\tau'\}\,$ and then we see that 
$\,1/|T_{s\wedge\tau'}|\,$ is a supermartingale, which implies that \parn 
\centerline{$2^{-1/2}\,\P (\tau'<\ii )\le \liminf_{s\to\ii}\limits\, \E ({1\over |T_{s\wedge\tau'}|}
\,1_{\{\tau'<\ii\}}) \le \liminf_{s\to\ii}\limits\, 
\E ({1\over |T_{s\wedge\tau'}|}) \le \E ({1\over |T_{0}|}) \le {1/2}\,$,} \parn 
and then that $\;\P (D<\ii ) \ge \P (\tau'=\ii ) \ge 1/\sqrt{2}\,$. \par 
   This concludes the proof of the second assertion in Theorem \ref{the.asymp}.\par\bigskip   

    3) {\it Existence of an asymptotic direction for the Schwarzschild diffusion, on $\,\{ D=\ii\}\,$.} 
\par\smallskip  

   We want to generalize the observation made in Section \ref{sec.relr} for $R=0$. Recall from Lemma
\ref{lem.cv} that it does not matter for this asymptotic behavior whether we consider the trajectories as
function of $\,s\,$ or of $t(s)\,$ (id est as viewed from a fixed point). \par   
   We shall use Section \ref{sec.relr} and Lemma \ref{lem.verif}, to proceed by comparison
between the flat Minkowski case $R=0$ and the Schwarzschild case $R>0$. \par 
   Let us split this proof into four parts. \par

\par \medskip  

 $(i)$ \ {\it  We have $\;{\ds\int_0^\ii {a_t\over r^2_t}\,dt <\ii}\;$ and $\;{\ds\int_0^\ii {U_t\over
r_t}\,dt <\ii}\;$, almost surely on $\,\{ D=\ii\}$..}\par\smallskip 

     We know from 1) above that $\,r_s\to\ii\,$ almost surely on $\,\{ D=\ii\}\,$. \par 

   The very beginning of this proof remains valid : Using (1,$i$) again, we have almost surely 
$\;{\ds\int_0^\ii {T^2_t\over a_t\,r^2_t}\,dt}\,$ and 
$\,{\ds\int_0^\ii (1-{\ts{3R\over 2r_t}})^2\, {b^2_t\over a_t\,r^3_t}\,dt}\;$ finite, whence 
$\;{\ds\int_0^\ii {b^2_t\over a_t\,r^3_t}\,dt}\;$ finite, and then, since 
$\; {\ds {a\over r^2}\le {T^2\over a\,r^2}+{b^2\over a\,r^4}+{1\over a\,r^2} }\;$, 
also $\;{\ds\int_0^\ii {a_t\over r^2_t}\,dt}\;$ finite, almost surely on $\,\{ D=\ii\}$. \par 
   Now by the unit pseudo-norm relation, we have 
$\;{\ds {U\over r}  = {b\over r^3} \le {a\over r^2\,\sqrt{1-{\ts{R\over r}}}}}\;$,  \ whence \parn 
$\;{\ds\int_0^\ii {U_t\over r_t}\,dt}\;$ finite, almost surely on $\,\{ D=\ii\}$.
\par\medskip 

$(ii)$ \ {\it  The perturbation of the Christoffel symbols due to $R$ is $\,\O(r\2)\,$.}
\par\smallskip 

   Recall from the beginning of Section \ref{sec.S} the values of the Christoffel symbols $\,\G^i_{jk}\,$. 
Denote by $\,\tilde\G^i_{jk}\,$ the difference between these symbols and their analogues for $R=0$, which
is a tensor, has only five non-vanishing components in spherical coordinates, and then is easily
computed in Euclidian coordinates
$\;(x_1=r\,\sin\f\,\cos\psi\,;\,x_2=r\,\sin\f\,\sin\psi\,;\,x_3=r\,\cos\f)\,$ : \ we find 
$$ \tilde\G^{x_i}_{x_j,x_k} = {\p x_i\over \p r}\times \Big( {\p r\over \p x_j}{\p r\over \p x_k}\G^r_{rr}
+ {\p \f\over \p x_j}{\p \f\over \p x_k}\G^r_{\f\f} + {\p \psi\over \p x_j}{\p \psi\over \p x_k}
\G^r_{\psi\psi}\Big) $$ 
$$ = {x_i\over r}\times \Big( {-R\over 2r(r-R)}{x_j\over r}{x_k\over r} + R{\p \f\over \p x_j} 
{\p \f\over\p x_k} + R\sin^2\f\,{\p \psi\over \p x_j}{\p \psi\over \p x_k}\Big) = \O(r\2)  $$ 
since $\;\Big|{\p\f\over\p x_j}\Big|\le 1/r\;$ and  $\;\Big|{\p\psi\over\p x_j}\Big|\le 1/(r\sin\f)\;$. 
The same is valid directly for the remaining components $\;\tilde\G^{x_i}_{t,t}\,$ and 
$\,\tilde\G^{t}_{x_j,t}\,$.  \par \medskip 

 $(iii)$ \ {\it The stochastic parallel transport converges, almost surely on $\,\{ D=\ii\}$.}
\par\smallskip 

   Let us denote by $\;\overleftarrow{\xi}(s)^i_j\;$ the matrix carrying out the inverse parallel
transport along the $C^1$ curve $\,(\xi_{s'}\,|\,0\le s'\le s)\,$, in the global pseudo-Euclidian
coordinates $\,(t,x_1,x_2,x_3)\,$. \par 
   We have 
$$ {d\over ds}\,\overleftarrow{\xi}(s)^i_j \, =\, \overleftarrow{\xi}(s)^i_k\times\G_{j\ell}^{k}(\xi_s)
\times\dot\xi_s^\ell\; , $$ 
so that, using $(ii)$ above and $\,|T_s|\le a_s\,$ : 
$$ \overleftarrow{\xi}(s)\, =\, \int_0^s \O\Big( r_v\2\times |\dot\xi_v|\Big)\,dv = 
\int_0^s \O\Big( |\dot t_v| + |\dot r_v| + r_v\,|\dot\t_v|\Big)\,r_v\2\,dv 
= \int_0^s\O\Big( 2\,\,{a_v\over r_v^2}+{U_v\over r_v}\Big)\, dv\, , $$ 
we can conclude by using $(i)$ above, that the Schwarzschild stochastic parallel transport (like its
inverse $\,\overleftarrow{\xi}(s)$) admits a finite limit as $s\to\ii\,$, almost surely on $\,\{ D=\ii\}$.
\par\medskip 

 $(iv)$ \ {\it End of the proof.} \par\smallskip 

    Let us consider $\;\eta_s := \overleftarrow{\xi}(s)\,\dot\xi_s\,$, \ for $\,s\ge 0\,$. \par 
This is a continuous process living on the fixed unit pseudo-sphere $\,T_{\xi_0}^1\M\,$. \parn  
Recall from Section \ref{sec.relg} that (for $\,0\le \ell\le 3\,$) 
$$ d\dot\xi_s^\ell = \s\sum_{k=1}^3 e^\ell_k(s)\circ dw^{k}_s
-\G_{ij}^\ell (x_s)\,\dot\xi_s^i\,\dot\xi_s^j\, ds  \; . $$ 
Therefore we get 
$$ d\eta_s^\ell\; =\; \s\sum_{k=1}^3 \overleftarrow{\xi}(s)^\ell_m\,e^m_k(s)\circ dw^{k}_s
- \overleftarrow{\xi}(s)^\ell_m\,\G_{ij}^m (x_s)\,\dot\xi_s^i\,\dot\xi_s^j\, ds  
+ (\circ d\overleftarrow{\xi}(s)^\ell_j)\,\xi_s^j  $$ 
$$ = \s\sum_{k=1}^3 \overleftarrow{\xi}(s)^\ell_m\,e^m_k(s)\circ dw^{k}_s
- \overleftarrow{\xi}(s)^\ell_m\,\G_{ij}^m (x_s)\,\dot\xi_s^i\,\dot\xi_s^j\, ds  
+ (\overleftarrow{\xi}(s)^\ell_m\,\G_{ji}^{m}(\xi_s)\,\dot\xi_s^i)\,\xi_s^j $$ 
$$ =\; \s\sum_{k=1}^3 \overleftarrow{\xi}(s)^\ell_m\,e^m_k(s)\circ dw^{k}_s 
=\; \s\sum_{k=1}^3 \tilde e^\ell_k(s)\circ dw^{k}_s \; , $$ 
where $\;\tilde e_k(s) := \overleftarrow{\xi}(s)\,e_k(s)\,$, for $\,1\le k\le 3\,$ and $\,s\ge 0\,$. \par  
   Observe that, for any $\,s\ge 0\,$, $\;\Big(\eta_s , \tilde e_1(s) , \tilde e_2(s) , \tilde
e_3(s)\Big)\;$ constitutes a pseudo-orthonormal basis of the fixed unit pseudo-sphere
$\,T_{\xi_0}^1\M\,$. \ Hence we find that the velocity process $\,\eta_\cdot\,$ defines a hyperbolic
Brownian motion on this unit pseudo-sphere $\,T_{\xi_0}^1\M\,$. \par 
   Now, according to Lemma \ref{lem.verif} and Section \ref{sec.relr}, we know that $\;\eta_s/a_s\;$
converges almost surely as $\,s\to\ii\,$ towards $\,(1,1,\tt_\ii)\,$ (in coordinates $(t,r,\t)$), for some
random $\,\tt_\ii\in\S^{2}$, exponentially fast. \par 
   Using $(iii)$ above, we deduce that $\;\dot\xi_s/a_s\;$, converges almost surely as $\,s\to\ii\,$
towards $\,(1,1,\hat\t_\ii)\,$, for some random $\,\hat\t_\ii\in\S^{2}$. \ 
This means also that the Schwarzschild diffusion seen from a fixed point, that is to say the implicit
trajectory $\;Z_\cdot := (t_s\mapsto (r_s,\t_s))\,$, sees almost surely its velocity $\,dZ_t/dt\,$
converging towards $\,(1,\hat\t_\ii)\,$, 1 being here the velocity of light. \par 
   Moreover, this shows a posteriori that $\,T_s\,$ goes to $+\ii$ and that $\,r_s\sim a_s\,$ as 
$\,s\to\ii\,$, and then that the convergences in $(i)$ and in $(iii)$ above occur exponentially fast, so
that it must be the same for the  convergences of $\;\dot\xi_s/a_s\;$ and $\,dZ_t/dt\,$. This allows to
integrate, to get finally the  generalization of Section \ref{sec.relr} to the Schwarzschild diffusion.

  This ends the whole proof of Theorem \ref{the.asymp}. \ $\;\diamond $ 
\bigskip \bigskip \bigskip 

\centerline{{\Large{\bf REFERENCES}}}  
\bigskip \bigskip 

\vbox{ \noindent 
{\bf [DF-C]} \ De Felice F. , Clarke C.J.S. \ \  {\it Relativity on curved manifolds. }
\par
\smallskip \hskip 11mm    Cambridge surveys on mathematical physics, Cambridge university press, 1990. }
\bigskip 

\vbox{ \noindent 
{\bf [D1]} \ Dudley R.M. \ \ {\it Lorentz-invariant Markov processes in relativistic phase space.}
\par \smallskip \hskip 33mm  Arkiv f\"or Matematik 6, n$^o$ 14, 241-268, 1965. }
\bigskip 

\vbox{ \noindent 
{\bf [D2]} \ Dudley R.M. \ \ {\it Asymptotics of some relativistic Markov processes.}
\par \smallskip \hskip 33mm  Proc. Nat. Acad. Sci. USA n$^o$ 70, 3551-3555, 1973. }
\bigskip 

\vbox{ \noindent 
{\bf [F-N]} \ Foster J. , Nightingale J.D. \ \  {\it A short course in General Relativity. }
\par
\smallskip \hskip 61mm    Longman, London 1979. }
\bigskip 

\vbox{  \noindent 
{\bf [I-W]} \ Ikeda N. ,  Watanabe S. \quad {\it Stochastic differential equations and diffusion
processes.}
\par \hskip 56mm North-Holland Kodansha, 1981. }
\medskip 

\vbox{ \noindent 
{\bf [S]} \ Stephani H. \ \ {\it General Relativity. } \qquad Cambridge university press, 1990. }
\bigskip 
\bigskip \bigskip 

\vbox{  
Jacques FRANCHI : \quad Universit\'e Louis Pasteur, I.R.M.A., 7 rue Ren\'e Descartes, \\ 
67084 Strasbourg cedex. France. \quad franchi@math.u-strasbg.fr 

\smallskip \noindent 

Yves LE JAN : \quad Universit\'e Paris Sud, Math\'ematiques, B\^atiment 425, 91405 Orsay. France. 
\quad  yves.lejan@math.u-psud.fr
}

\end{document}